\documentclass[a4paper]{article}
\usepackage{amsmath,amssymb}
\usepackage[frenchb]{babel}
\usepackage[utf8]{inputenc}
\usepackage[dvips]{epsfig}

\newcommand{\QQ}{\mathbb{Q}}
\newcommand{\CC}{\mathbb{C}}

\newcommand{\qzeta}{\boldsymbol{\zeta}_q}
\newcommand{\qL}{\mathbf{L}_q}

\newtheorem{theorem}{Théorème}[section] 
\newtheorem{proposition}[theorem]{Proposition}

\newtheorem{lemma}[theorem]{Lemme}

\title{Fractions de Bernoulli-Carlitz \\ et opérateurs q-Zeta}
\author{F. Chapoton}
\date{\today}

\begin{document}

\maketitle

\begin{abstract}
  On introduit une déformation des séries de Dirichlet d'une variable
  complexe $s$, sous la forme d'un opérateur pour chaque nombre
  complexe $s$, agissant sur les séries formelles en une variable $q$
  sans terme constant. On montre que les fractions de
  Bernoulli-Carlitz sont les images de certains polynômes en $q$ par
  les opérateurs associés à la fonction $\zeta$ de Riemann aux entiers
  négatifs.
\end{abstract}

\section{Opérateurs associés aux séries de Dirichlet}

A chaque suite $(a_{n})_{n\geq 1}$ de nombres complexes, on associe la
série de Dirichlet formelle
\begin{equation}
  \zeta(a,s) = \sum_{n \geq 1} \frac{a_n}{n^s},
\end{equation}
où $s$ est un nombre complexe. Le produit de deux telles séries
correspond à la convolution des suites de nombres complexes :
\begin{equation}
  \zeta(a,s)\zeta(b,s)=\zeta(a*b,s),
\end{equation}
où $(a*b)_n=\sum _{d|n} a_d \, b_{n/d}$.

Considérons maintenant l'anneau $\CC[[q]]$ des séries formelles à
coefficients complexes en une variable $q$. Pour $n\geq 1$, on définit le nombre quantique 
\begin{equation}
  [n]=\frac{q^n-1}{q-1}
\end{equation}
et l'opérateur de Frobenius formel 
\begin{equation}
  F_n(f)=f(q^n),
\end{equation}
agissant sur $\CC[[q]]$ et par restriction sur le sous espace
$q\,\CC[[q]]$ des séries formelles sans terme constant.

A chaque suite $(a_n)_{n\geq 1}$, on peut alors associer l'opérateur
\begin{equation}
  \qzeta(a,s) = \sum_{n \geq 1} \frac{a_n F_n}{[n]^s},
\end{equation}
qui agit sur $q\,\CC[[q]]$. Cet opérateur est bien défini pour tout
nombre complexe $s$. 

\begin{lemma}
  \label{commute}
  Pour tous $m,n$ entiers non nuls, on a la relation 
  \begin{equation}
    \frac{F_m}{[m]^s}\frac{F_n}{[n]^s}=\frac{F_{m n}}{[m n]^s}.
  \end{equation}
  En particulier, ces deux opérateurs commutent.
\end{lemma}

Par conséquent, 
\begin{proposition}
  L'application $a \mapsto \qzeta(a,s)$ est un morphisme de l'algèbre
  des suites de nombres complexes pour la convolution dans l'algèbre
  des opérateurs linéaires sur $q\,\CC[[q]]$.
\end{proposition}

En particulier, si la suite $a$ est multiplicative, l'opérateur
$\qzeta(a,s)$ admet un produit eulérien, exactement similaire au
produit eulérien de la série de Dirichlet $\zeta(a,s)$.

Par exemple, la fonction $\zeta$ de Riemann est la série de Dirichlet
associée à la suite constante égale à $1$. On a donc un opérateur
$\qzeta$ défini par
\begin{equation}
  \label{def_zetaq}
    \qzeta(s) = \sum_{n \geq 1} \frac{F_n}{[n]^s}.
\end{equation}
Cet opérateur admet le produit eulérien
\begin{equation}
     \qzeta(s) = \prod_{p} \left( 1-\frac{F_p}{[p]^s} \right)^{-1},
\end{equation}
où le produit porte sur l'ensemble des nombres premiers.

\textbf{Remarque :} un opérateur Zeta proche de l'opérateur $\qzeta(0)$ a été
introduit dans \cite{meyer}. Il utilise plutôt la variable $\tau$
telle que $q=\exp(2 i \pi \tau)$.
 
\section{Fractions de Bernoulli-Carlitz}

Introduisons maintenant les fractions de Bernoulli-Carlitz, qui vont
jouer le rôle des nombres de Bernoulli. Ces fractions ont été définies
par Carlitz \cite{carlitz1,carlitz2}, puis étudiées par Koblitz et
Satoh \cite{koblitz,satoh}. Elles sont aussi apparues plus récemment
dans \cite{qomega}, avec une motivation algébrique.

On peut les définir par la récurrence suivante : on pose $\beta_0=1$ et
\begin{equation}
  \label{recu_beta}
  q(q \beta +1)^n-\beta_n=
  \begin{cases}
    1 \text{ si }n=1,\\
    0 \text{ si }n>1.
   \end{cases}
\end{equation}
Dans cette expression, par convention, on pose $\beta^i=\beta_i$ pour
tout entier $i$, après avoir développé la puissance du binôme.

On obtient ainsi une suite de fractions rationnelles $\beta_n$ dans le
corps $\QQ(q)$. On peut montrer (voir \cite{carlitz1}) que ces
fractions n'ont pas de pôle en $q=1$ et que la valeur de $\beta_n$ en $q=1$
est le nombre de Bernoulli $B_n$.

Les premières fractions sont
\begin{align}
  \beta_0&=1,\\
  \beta_1&=-\frac{1}{\Phi_2},\\
  \beta_2&=\frac{q}{\Phi_2 \Phi_3},\\
  \beta_3&=\frac{q(1-q)}{\Phi_2\Phi_3\Phi_4},\\
  \beta_4&=\frac{q(q^4 - q^3 - 2 q^2 - q + 1)}{\Phi_2\Phi_3\Phi_4\Phi_5},
\end{align}
où les $\Phi_n$ sont les polynômes cyclotomiques.

Soit $\mathbb{B}(t)$ la série génératrice des $\beta_n$ :
\begin{equation}
  \label{serie_gen}
  \mathbb{B}(t)=\sum_{n\geq 0} \beta_n \frac{t^n}{n!}.
\end{equation}
La récurrence \eqref{recu_beta} se traduit en une équation fonctionnelle pour $\mathbb{B}$ :
\begin{equation}
  \label{eq_fonct}
\mathbb{B}(t)=q e^t \mathbb{B}(qt)+1-q-t.  
\end{equation}

On vérifie (en utilisant la relation $1 +q [n]=[n+1]$) que la solution
de \eqref{eq_fonct} est donnée par
\begin{equation}
\mathbb{B}(t)=\sum_{n\geq 0} (1-q) q^n e^{[n]t} - t q^{2 n} e^{[n]t}.
\end{equation}

Par dérivations successives par rapport à $t$, on obtient la formule explicite
\begin{equation}
  \beta_n=\sum_{k \geq 1}q^k [k]^{n-1}-(n+1)\sum_{k \geq 1}q^{2 k} [k]^{n-1},
\end{equation}
valable pour $n\geq 2$.

En utilisant l'opérateur $\qzeta(1-n)$, on obtient donc le résultat suivant.

\begin{theorem}
  Pour $n\geq 2$, on a
  \begin{equation}
    \label{q-relation}
    \beta_n=\qzeta(1-n) \left( q-(n+1) q^2 \right).
  \end{equation}
\end{theorem}

L'égalité \eqref{q-relation} est un $q$-analogue naturel du résultat d'Euler qui relie les valeurs prises aux entiers négatifs par la fonction $\zeta$ de Riemann et les nombres de Bernoulli $B_n$ : pour tout $n \geq 2$, on a
\begin{equation}
  B_n=\zeta(1-n) (-n).
\end{equation}

On peut espérer déduire le résultat classique d'Euler par un passage à
la limite $q=1$ en un sens approprié. Un travail analytique reste à
faire pour obtenir cela.



\section{Opérateur $q$-Zeta de Hurwitz}

Comme dans le cas classique, les résultats précédents ont des
analogues pour l'opérateur $q$-Zeta de Hurwitz.

Soit $x \in \QQ_{>0}$. On note encore $F_{n+x}$ l'opérateur de Frobenius
formel
\begin{equation}
  F_{n+x}(f)=f(q^{n+x}),
\end{equation}
qui agit sur les séries de Puiseux sans terme constant en la variable $q$. 

On définit alors l'opérateur Zeta de Hurwitz
\begin{equation}
  \label{def_hurwitz}
    \qzeta(s,x) = \sum_{n \geq 0} \frac{F_{n+x}}{[n+x]^s}.
\end{equation}

\begin{proposition}
  Pour tout $N \geq 1$ et tout $s$ dans $\CC$, on a la relation de distribution
  \begin{equation}
    \sum_{0\leq j < N} \frac{F_N}{[N]^s}\qzeta(s,\frac{x+j}{N})=\qzeta(s,x).
  \end{equation}
\end{proposition}

La preuve est immédiate, en utilisant la définition
\eqref{def_hurwitz} et une version adaptée du Lemme \ref{commute}.

Par ailleurs, on introduit (d'après Carlitz) des $q$-polynômes de
Bernoulli par la formule symbolique
\begin{equation}
  \beta_n(x)=(q^x \beta+[x])^n.
\end{equation}
Ce sont des polynômes en la variable $q^x$ à coefficients dans $\QQ(q)$.

La série génératrice des $q$-polynômes de Bernoulli est alors
\begin{equation}
  e^{[x]t}\mathbb{B}(q^x t).
\end{equation}

La série génératrice des $\beta_n(x)$ est donc explicitement donnée par
\begin{equation}
e^{[x]t} \sum_{n\geq 0} (1-q) q^n e^{[n]q^x t} - t q^{2 n+x} e^{[n]q^x t}.
\end{equation}
soit encore (par la relation $[x]+q^x[n]=[n+x]$)
\begin{equation}
\sum_{n\geq 0} (1-q) q^n e^{[n+x] t} - t q^{2 n+x} e^{[n+x] t}.
\end{equation}

Par dérivations successives par rapport à $t$, on obtient la formule explicite
\begin{equation}
  q^x \beta_n(x)=\sum_{k \geq 0}q^{k+x} [k+x]^{n-1}-(n+1)\sum_{k \geq 0}q^{2 k+2 x} [k+x]^{n-1},
\end{equation}
valable pour $n\geq 0$.

En exprimant ceci avec l'opérateur $q$-Zeta de Hurwitz, on trouve l'égalité suivante.

\begin{proposition}
  \label{relation_hurwitz}
  Pour $n \geq 0$, on a la relation
  \begin{equation}
    q^x \beta_{n}(x)=\qzeta(1-n,x) \left( q-(n+1) q^2 \right).
  \end{equation}
\end{proposition}

\section{Opérateurs $q$-$L$ de Dirichlet}

Comme dans le cas classique, on peut généraliser les résultats
précédents aux opérateurs $q$-$L$ associés aux caractères de
Dirichlet, en utilisant l'opérateur $q$-Zeta de Hurwitz.

Soit $\chi$ un caractère de Dirichlet de conducteur $N$. On définit un opérateur
\begin{equation}
    \qL(\chi,s) = \sum_{n \geq 1} \frac{\chi(n) F_n}{[n]^s}.
\end{equation}

On suppose maintenant que $\chi$ n'est pas un caractère trivial, pour
simplifier.

\begin{proposition}
  \label{zeta_zeta}
  On a
  \begin{equation}
    \qL(\chi,s)=\sum_{0\leq j <N}\chi(j) \frac{F_N}{[N]^{s}} \qzeta(s,j/N).
  \end{equation}
\end{proposition}

La preuve est similaire à celle de la relation de distribution.

On introduit alors les analogues des fractions de Bernoulli-Carlitz
pour le caractère $\chi$ :
\begin{equation}
  \beta_{\chi,n}= \sum_{0\leq j< N} \chi(j) \frac{F_N}{[N]^{1-n}}\left (q^{j/N} \beta_n(j/N) \right).
\end{equation}

On déduit des propositions \ref{relation_hurwitz} et \ref{zeta_zeta}
l'égalité suivante.

\begin{proposition}
  Pour $n\geq 0$, on a
  \begin{equation}
    \label{chi-relation}
    \beta_{\chi,n}=\qL(\chi,1-n) \left( q-(n+1) q^2 \right).
  \end{equation}
\end{proposition}

\section{Calcul des valeurs aux entiers négatifs}

Soit $\chi$ un caractère de Dirichlet et $r$ un entier strictement
positif. Comme $\chi$ est périodique, la somme $f_\chi(r)=\sum_{n\geq
  1} \chi(n) q^{n r}$ est une fraction rationnelle. C'est aussi la
valeur $\qL(\chi,0) q^r$. Par exemple, pour le caractère trivial,
on obtient $q^r/(1-q^r)$.

Si $f$ est une fonction de la variable $r$, on note $\Delta$
l'opérateur suivant (qui est une forme de $q$-différence) :
\begin{equation}
  \Delta(f)=\frac{f(r+1)-f(r)}{q-1}.
\end{equation}

L'opérateur $\Delta$ agit sur l'ensemble des fractions en $q$ et $q^r$
qui, en tant que fonction de $q$ pour tout $r$ entier fixé, ont des
pôles seulement sur le cercle unité et sont nulles en zéro.

L'opérateur $\Delta$ vérifie
\begin{equation}
  \Delta(q^{nr})=[n] q^{nr}.
\end{equation}

On a donc, pour tout $i\leq 0$ et tout entier $r\geq 1$, la valeur
\begin{equation}
  \qL(\chi,i) q^r=\Delta^i f_{\chi}.
\end{equation}

Il en résulte donc que toutes ces valeurs sont des fractions
rationnelles en $q$, nulles en zéro et avec pôles sur le cercle unité.


\begin{proposition}
  Soit $P$ un polynôme en $q$ sans terme constant. Si $i\leq 0$, alors
  $\qzeta(\chi,i)P$ est une fraction rationnelle en $q$.
\end{proposition}

\section{Remarques et spéculations}

Parmi les différentes fonctions $q$-Zeta qui sont apparues dans la
littérature \cite{kaneko,tsumura,kim}, certaines sont obtenues par
application de l'opérateur $q$-Zeta à diverses puissances de $q$. Il
paraît plus naturel d'étudier l'opérateur lui-même plutôt qu'un choix
arbitraire de ses valeurs. 

On peut imaginer que les valeurs spéciales des fonctions $L$
correspondent à l'évaluation des opérateurs $q$-$L$ en des arguments
particuliers. Il faudrait écrire la valeur spéciale comme un
quotient et trouver des $q$-analogues du numérateur et du
dénominateur.

Vu son produit eulérien, on peut aussi imaginer l'existence d'une équation
fonctionnelle pour l'opérateur $\qzeta(s)$. Il s'agirait de compléter
le produit eulérien en introduisant un opérateur $\boldsymbol{\gamma}(s)$
correspondant au facteur archimédien. Cet opérateur devrait commuter
avec tous les opérateurs $\frac{F_n}{[n]^s}$, pour que le produit
eulérien complété reste totalement commutatif.

\medskip

On observe numériquement la propriété suivante des fractions de
Bernoulli-Carlitz. Leur numérateur a quelques zéros réels positifs,
beaucoup de zéros sur le cercle unité et quelque paires de zéros
complexes. Ceci semble aussi vrai pour les fractions similaires
associées aux caractères de Dirichlet.



\bibliographystyle{alpha}
\bibliography{qzeta}

\end{document}